\documentstyle[mathrsfs, times, a4paper]{amsart}

\theoremstyle{plain}
\newtheorem{Thm}{Theorem}[subsection]
\newtheorem*{MTH}{Main Theorem}

\newtheorem*{Kor}{Corollary}
\newtheorem{Prop}[Thm]{Proposition}
\newtheorem{Lem}[Thm]{Lemma}
\theoremstyle{remark}
\newtheorem{Rem}[Thm]{Remark}
\newtheorem{Ex}[Thm]{Example}
\newtheorem*{Loc}{Local definition of no future importance}
\newcommand{\Amp}{\mathop{\rm Amp}}
\newcommand{\Num}{\mathop{\rm Num}}
\newcommand{\chh}{{\mathop{\rm ch}}}
\newcommand{\ch}{{\mathop{\rm hilb}}}
\newcommand{\td}{{\mathop{\rm td}}}
\newcommand{\eps}{\varepsilon}
\newcommand{\catqot}{/\hskip-3pt/}
\newcommand{\C}{{\Bbb C}}
\newcommand{\E}{{\mathscr E}}
\newcommand{\Ext}{\mathop{\rm Ext}}

\newcommand{\F}{{\mathscr F}}
\newcommand{\G}{{\mathscr G}}
\newcommand{\Hom}{\mathop{\rm Hom}}

\newcommand{\id}{\mathop{\rm id}}
\renewcommand{\L}{{\mathscr L}}

\renewcommand{\O}{{\mathscr O}}
\renewcommand{\P}{{\Bbb P}}
\newcommand{\Pic}{\mathop{\rm Pic}}
\newcommand{\Q}{{\Bbb Q}}
\newcommand{\R}{{\Bbb R}}
\newcommand{\SL}{\mathop{\rm SL}}
\newcommand{\GL}{\mathop{\rm GL}}
\newcommand{\Z}{{\Bbb Z}}
\newcommand{\N}{{\Bbb N}}
\newcommand{\la}{\lambda}
\newcommand{\lra}{\longrightarrow}

\newcommand{\p}{\prime}
\newcommand{\q}{\quad}
\renewcommand{\phi}{\varphi}
\newcommand{\rk}{\mathop{\rm rk}}

\title{Walls for Gieseker semistability and the Mumford-Thaddeus
principle\\ for moduli spaces of sheaves\\
over higher dimensional bases}
\thanks{Currently supported by a grant of the Emmy Noether Institute
at Bar-Ilan University}

\author{Alexander Schmitt}

\begin{document}
%
\abstract{Let $X$ be a projective manifold over $\C$. Fix two
ample line bundles $H_0$ and $H_1$ on $X$. It is the aim of this note
to study the variation of the moduli spaces of Gieseker semistable
sheaves for polarizations lieing in the cone spanned by $H_0$ and $H_1$.
We attempt a new definition of walls which naturally describes the
behaviour of Gieseker semistability. By means of an example, we
establish the possibility of non-rational walls which is a
substantially new phenomenon compared to the surface case. Using the
approach by Ellingsrud and G\"ottsche via parabolic sheaves, we were
able to show that the moduli spaces undergo a sequence of GIT flips
while passing a rational wall. We hope that our results will be
helpful in the study of the birational geometry of moduli spaces
over higher dimensional bases.}
\endabstract

\maketitle

\markboth{Alexander Schmitt}{MUMFORD-THADDEUS PRINCIPLE}

\section*{Introduction}
Fix an $n$-dimensional
smooth projective manifold $X$ over the complex
numbers as well as a function $p\colon \Num(X)\lra\Z$,
called \it Hilbert form\rm .
Define
$N^1_{\Q}(X):=\Num(X)\otimes_{\Z}\Q$ and similarly
$N^1_{\R}(X)$, and finally let $\Amp_{\Q}(X)$ and $\Amp_{\R}(X)$
be the cones in $N^1_{\Q}(X)$ and $N^1_{\R}(X)$, resp., spanned
by the classes of ample line bundles.
Assuming that $H$ is the class of an ample line bundle,
we define $P_H(\E)$ as the polynomial such that $P_H(\E)(n)=
\chi(\E\otimes H^{\otimes n})$ for any natural number $n$.
The sheaf $\E$ is then called \it Gieseker $H$-(semi)stable \rm
(or just {\it $H$-(semi)stable})
if and only if every non-zero proper subsheaf $\F$ of $\E$ satisfies
$P_H(\F)/\rk\F\q(\le)\q P_H(\E)/\rk\E$. There is a projective moduli
space
${\mathscr M}_H:={\mathscr M}_H(p)$ of S-equivalence classes of
Gieseker $H$-semistable torsion free coherent sheaves $\E$ with
Hilbert form $p$, i.e., $p([D])=\chi(\E\otimes \O_X(D))$  for all
$[D]\in\Num(X)$. Note that this determines the rank of $\E$,
henceforth denoted by $r$, the numerical equivalence class of $c_1\E$,
henceforth denoted by $c_1$, and $c_2\E$ as a linear form
on the subvectorspace of $H^{2n-4}(X,\Q)$ spanned by $(n-2)$-fold
intersections of divisors, as such it is called $c_2$.
By its very definition, the space ${\mathscr M}_H$ depends on the chosen
polarization, and it is an interesting and important problem to compare
${\mathscr M}_{H_0}$ to  ${\mathscr M}_{H_1}$ for different
polarizations $H_0$
and
$H_1\in \Amp_{\Q}(X)$.
For surfaces, this problem has been
thoroughly studied. A brief discussion of this topic and appropriate
references can be found in \cite{HL}. The most general result
in this direction has been obtained in \cite{MW} where it is shown
that the moduli spaces are related by a sequence of GIT flips.
A similar result can be obtained using moduli spaces of parabolic
sheaves as mentioned in the paper \cite{EG}. In this note we aim
at a generalization of the results of \cite{MW} to higher dimensions,
using the approach of \cite{EG}.
However, there arise new problems due to the appearence of walls
which do not lie in $N^1_{\Q}(X)$. Our result is summarized in the
following
\begin{MTH}
Given two polarizations $H_0$ and $H_1$, there is a finite subset
$w$ of $\Delta:=\bigl\{\, (1-\la)H_0+\la H_1\ |\ \la\in[0,1]\,\bigr\}$
such that the notion of Gieseker (semi)stability remains constant
within each connected component of $\Delta\setminus w$.
If the polarization passes through a wall of $w\cap N_\Q^1(X)$, then
the moduli spaces undergo a sequence of $\C^*$-flips.
\end{MTH}
In the case of crossing a real wall,
one cannot expect such a result, because it would yield an algebro
geometric construction of a moduli space of Gieseker semistable sheaves
w.r.t.~a \sl real \rm polarization which seems most unlikely in my eyes.
However, in this case,
some suitable fibre spaces over the moduli spaces can be obtained by
a sequence of $\C^*$-flips from the same Quot scheme. This will be
explained in Section~\ref{OST}.
\par
In general, the hope is that ${\mathscr M}_{H_0}$ and ${\mathscr M}_{H_1}$
will be --- under suitable assumptions --- birational to each other,
although other results
indicate that moduli spaces over higher dimensional bases are not at all
well-behaved, e.g., they can have arbitrarily many components
(\cite{E},\cite{BM}).
The flips between the moduli spaces can be very helpful in this
context.
In fact, one should be able to obtain quite explicit descriptions
of the exceptional sets of the flips. Then, one is left
with estimating the dimension of these exceptional sets,
and this might
be the hard part.
\par
In the case of crossing a rational wall,
our construction gives the following: There is a quasi-projective
scheme ${\frak X}$, an ample line bundle ${\frak L}$ on ${\frak X}$,
and a $\C^*$-action on ${\frak X}$ together with two linearizations
$\sigma_0$ and $\sigma_1$ of this action in ${\frak L}$ such that
${\frak X}\catqot_{\sigma_{0,1}}\C^*= {\mathscr M}_{H_{0,1}}$.
Let ${\frak X}_i$, $i=1,...,t$, be the irreducible components
of ${\frak X}$. Since $\C^*$ is irreducible, the action 
preserves those components.
So, the ${\mathscr M}_{0,1}^i$ will be the irreducible components of 
${\mathscr M}_{H_{0,1}}$, $i=1,...,t$. By general properties of $\C^*$-actions
(e.g.\ \cite{Th2}, \cite{OST}) one gets
\begin{Kor} 
Under the above hypotheses, if for $i_0\in \{\, 1,...,t\,\}$
both ${\mathscr M}_0^{i_0}$ and ${\mathscr M}_1^{i_0}$ are non-empty,
then they are birationally equivalent.
\end{Kor} 
\subsection*{Acknowledgements}
The paper was inspired by the suggestion of Professor Mir\'o-Roig
to study the birational geometry of moduli spaces of sheaves over
higher dimensional base varieties. During the preparation of it,
I profited from discussions with Laura Costa and Manfred Lehn.
During the preparation of the article, 
the author was supported by grant \#1996SGR00060 of the Direcci\'o General
de Recerca, Generalitat de Catalunya.
%
\section{Preparations}
\subsection{Walls for slope semistability}
For technical reasons, we will have to consider the notion
of slope semistability for all $H\in\Amp_{\R}(X)$.
So, let us fix such an $H$. For a torsion free coherent sheaf
$\E$, define its {\it $H$-slope} as $\mu_H\E:=c_1\E.H^{n-1}/\rk\E$
and call $\E$ \it slope $H$-(semi)stable \rm if $\mu_H\F\ (\le)\
\mu_H\E$ for any non-zero proper subsheaf $\F$ of $\E$.
\begin{Ex}
\label{zb1}
Let $X\subset \P_2\times\P_2$ be a smooth hypersurface
in $|\O(1,1)|$. The nef cone of $X$ is spanned by
$H_0:=\pi_1^*\O_{\P_2}(1)$ and $H_1:=\pi_2^*\O_{\P_2}(1)$.
Set $H_\la=(1-\la)H_0+\la H_1$.
We have $H_0^3=0=H^3_1$ and $H_0^2.H_1=1=H_0.H_1^2$.
Define $E:=\O(2,-1)\oplus \O(-2,1)$. This bundle will be
slope $H_\la$-semistable if and only if
$$
0\q=\q\O(2,-1).H_{\la}^2\q=\q -\la^2+4\la-1.
$$
This equation has the  (irrational) solutions $\la_{\pm}:=2\pm
\sqrt{3}$. Note that $\la_-$ gives a real class in the ample cone.
Hence, $E$ is semistable only with respect to a single real class!
Thus, the study of sheaves which are slope semistable
w.r.t.~a real
class cannot necessesarily be reduced to the study of vector bundles
which are slope semistable for some rational class.
\end{Ex}
Fix two polarizations $H_0$ and $H_1$ in $\Amp_{\Q}(X)$ and denote
the line segment joining them by $\Delta$.
In this section, $H_\la$ stands for the
polarization $(1-\la) H_0 + \la H_1$,
$\la\in [0,1]$.
We are interested
in the family ${\frak F}(\Delta)$ of isomorphy classes of torsion free
coherent sheaves $\E$ with Hilbert form $p$
for which there exists a rational polarization $H\in\Delta\cap
N^1_{\Q}(X)$ w.r.t.~which $\E$ is slope semistable.
\par
For any sheaf $\E$
and any non-zero
proper subsheaf
$
\F\subset\E$ define $\xi_{\F,\E}:=[c_1\F/\rk\F-c_1/r]$.
We begin with the following observation.
\begin{Lem}
\label{Obs1}
Let $\la_0\in [0,1)$ and $\la_1\in (\la_0,1)$. Denote
the family of isomorphy classes of slope $H_{\la_0}$-semistable
torsion
free coherent
sheaves with Hilbert form $p$ by ${\frak F}(H_{\la_0})$.
Then there is a constant $C$ such that for any $\E$ with $[\E]\in
{\frak F}(H_{\la_0})$ and any non-zero proper subsheaf $\F$ of $\E$
the condition $\xi_{\F,\E}.H_{\la_0}^{n-1}< C$ implies
$\xi_{\F,\E}.H_{\la}^{n-1}< 0$ for all $\la\in[\la_0,\la_1]$.
\end{Lem}
\begin{pf}
We may assume $\la_0=0$.
Then ${H}_\la^{n-1}=\sum_{i=0}^{n-1}
{n-1\choose i} (1-\la)^i
\la^{n-i-1}
H_{0}^iH_{1}^{n-1-i}$.
Since ${\frak F}(H_0)$ is a bounded family, there are constants
$K_0,...,K_{n-2}$ such that
$\xi_{\F,\E}.H_{0}^i\allowbreak H_{1}^{n-1-i}\le K_i$, $i=0,...,n-2$,
for all $\E$ with $[\E]\in{\frak F}(H_0)$ and all subsheaves
$0\neq
\F\subset\E$.
Setting $K:=\max\bigl\{\, \sum_{i=0}^{n-2} {n-1\choose i}
(1-\la)^i
\la^{n-i-1}
K_i\ |\ \la\in [0,1]\,\bigr\}$, we conclude
that,
for $\la\in [0,\la_{1}]$,
$\E$ with $[\E]\in{\frak F}(H_0)$, and all subsheaves
$0\neq
\F\subset\E$,
$$
0\ \le\  \xi_{\F,\E}.{H}_\la^{n-1}\ \le\
(1-\la)^{n-1}
\xi_{\F,\E}.{H}_{0}^{n-1}+K
\ \le\
(1-\la_1)^{n-1}
\xi_{\F,\E}.{H}_{0}^{n-1}+K
$$
implies $\xi_{\F,\E}.{H}_{0}^{n-1}\ge
-K/(1-\la_1)^{n-1}$, and we are done.
\end{pf}
As important consequence, we note
\begin{Prop}
\label{wichtig1}
Let $\E$ be a torsion free coherent sheaf such that
$[\E]\in{\frak F}(\Delta)$.
Suppose that $\E$ is slope semistable w.r.t.~$H_{\la_0}$ with
$\la_0\neq 1$.
Assume that for any subsheaf $\F\subset \E$
there is an open neighborhood $U\subset
[0,1]$ of $\la_0$, such that $\xi_{\F,\E}.H_{\la}^{n-1}\le 0$ for all
$\la\in U$. Then either $\E$ is slope $H_\la$-semistable for any
$\la\in[\la_0,1]$, or there exists a
number
$\la_+>\la_0$ such that {\rm i)} $\E$ is slope $H_\la$-semistable for
any
$\la\in
[\la_0,\la_+]$, {\rm ii)} there exists a  saturated
non-zero proper
subsheaf
$\F_+\subset\E$ with $\mu_{H_{\la_+}}\F_+=\mu_{H_{\la_+}}\E$
such that
$
\bigl( (\rk\F_+-1) c_1^2\F_+ - 2\rk\F_+ c_2\F_+\bigr).H^{n-2}_{\la_+}
\q \le
\q 0,
$
and, for $\G_+:=\E/\F_+$,
$
\bigl( (\rk\G_+-1) c_1^2\G_+ - 2\rk\G_+ c_2\G_+\bigr).H^{n-2}_{\la_+}
\q\le
\q 0,
$
{\rm iii)} $\E$ is
not slope $H_{\la}$-semistable for
$\la>\la_+$ close enough.
\end{Prop}
\begin{Rem} i)
Likewise, one can construct under the assumption $\la_0\neq 0$
a number $\la_-<\la_0$ and a subsheaf $\F_-$ with the respective
properties.
\par
ii) The sheaf $\F_+$ is not necessarily slope $H_\la$-desemistabilizing
for $\la>\la_+$.
\par
iii) The need for this proposition arises from the fact that I don't
know if the Bogomolov inequality continues to hold for real
polarizations.
\end{Rem}
\begin{pf}
We may suppose that $\E$ is not slope $H_{\la_1}$-semistable for
some rational $\la_1>\la_0$.
If a subsheaf $\F$ slope desemistabilizes
$\E$ for some $H_\la$ with $\la\in [\la_0,\la_1]$,
then we must have $\xi_{\F,\E}.H_{\la_0}^{n-1}\ge C$,
by Lemma~\ref{Obs1}.
The set of saturated subsheaves $\F$ of $\E$
with
$\xi_{\F,\E}.H_{\la_0}^{n-1}\ge C$ is bounded (\cite{HL}, Lem.~1.7.9).
In particular, there are only finitely many elements $\xi$
in $(1/r!)\Num(X)$ of the form $\xi_{\F,\E}$ for which there
is a $\la\in[\la_0,\la_1]$ with
$\xi.H_{\la}^{n-1}\ge 0$.
Denote these elements by $\xi_1,...,\xi_\nu$
and set $f_i(\la):= \xi_i.H_\la^{n-1}$. Let $\la_+$ be the smallest
number in $(\la_0,\la_1]$ at which one of the polynomial
functions $f_i(\la)$ undergoes a change of sign.
Then, by construction, $\E$ is slope semistable for all
$H_{\la}$ with $\la\in[\la_0,\la_+)$, properly slope
$H_{\la_+}$-semistable, and slope unstable for values $\la>\la_+$,
close enough.
Furthermore, without loss of generality, we can
assume that $f_1(\la)\ge f_i(\la)$ for $i=2,...,\nu$ and all
$\la\le \la_+$, close enough to $\la_+$.
Pick some saturated subsheaf $\F_+$ such that
$f_1(\la)=\xi_{\F_+,\E}.H^{n-1}_\la$. Then $\F_+$ and the quotient
$\G_+$ are by construction slope $H_{\la}$-semistable for $\la\le\la_+$,
close enough, and $\mu_{H_{\la_+}}\F_+=\mu_{H_{\la_+}}\E
=\mu_{H_{\la_+}}\G_+$. Furthermore, $\F_+$ and $\G_+$
satisfy
$$
\bigl( (\rk\F_+-1) c_1^2\F_+ - 2\rk\F_+ c_2\F_+\bigr).H^{n-2}_{\la}
\q \le
\q 0
$$
and
$$
\bigl( (\rk\G_+-1) c_1^2\G_+ - 2\rk\G_+
c_2\G_+\bigr).H_\la^{n-2}
\q\le
\q 0
$$
for all rational $\la\le\la_+$, close enough, by the Bogomolov theorem
(\cite{HL}, Thm.~7.3.1). Thus, the proposition is proved.
\end{pf}
\begin{Ex}
\label{zb2}
This time, we consider a smooth hypersurface $X\subset\P_2\times\P_2$
in the linear system $|\O(3,3)|$. Using notations analogous to
those in Example~\ref{zb1}, we have generators $H_0$ and $H_1$
of the nef cone of $X$ with $H_0^3=0=H_1^3$ and $H_0^2.H_1=3=H_0.H_1^2$.
The space $X$ is a Calabi-Yau threefold with
$c_2(X)=3H_0^2+3H_1^2+9
H_0H_1$. First, we check that there is a non-split extension
$$
0\lra \O_X(3,0)\lra E\lra \O_X(0,1)\lra 0.
$$
Such extensions are parametrized by  the 
$\Ext^1(\O_X(0,1),\O_X(3,0))=\allowbreak H^1(\O_X(3,-1))$.
Observe $h^0(\O_X(3,-1))=0=h^0(\O_X(-3,1))=h^3(\O_X(3,-1))$, so that
Riemann-Roch gives
$-h^1(\O(3,-1))\le (1/6)(3H_0-H_1)^3+(1/12)(3H_0-H_1).
(3H_0^2+3H_1^2+9 H_0H_1)=-3$.
Besides subsheaves of $\O(3,0)$, $E$ could have subsheaves of the
form $\O(-k,1)$ with $k\ge 1$, because the extension does not split.
Subsheaves of the latter form do not destabilize if
$\xi_{\O(-1,1),E}.H_\la^2<0$ where
$\xi_{\O(-1,1),E}=-(5/2) H_0 + (1/2) H_1$.
One checks that this is fulfilled for all $\la> \la^*:=(5/4)-(\sqrt{21}/4)$.
Thus, for $\la>\la^*$, the middle term $E$ of such a non-split extension
is slope $H_{\la}$-(semi)stable if and only if $\O_X(3,0)$ does
not de(semi)stabilize $E$.
We have $\xi:=\xi_{\O_X(3,0),E}= (3/2)H_0-(1/2)H_1$, and the equation
$\xi.H_\la^{2} (\le)0$ reads
$$
{3\over 2}\bigl( -2\la^2+6\la-1\bigr)\qquad (\le)\qquad 0.
$$
Thus, $E$ is slope stable for all polarizations $H_{\la}$
with $\la^*<\la<(3/2)-(1/2)\sqrt{7}$, properly slope semistable for
$H_{(3/2)-(1/2)\sqrt{7}}$,
and not semistable for any polarization $H_\la$ with 
$\la> (3/2)-(1/2)\sqrt{7}$.
\end{Ex}
\begin{Rem}
This example exhibits an interesting phenomenon.
Although our set-up is completely algebro-geometric, we naturally
encounter objects which are not readily accessible by algebraic
methods. In particular, it becomes clear that in order to completely
solve our problem we have to find the right notion of Gieseker
semistability w.r.t.~an arbitrary K\"ahler class and to construct
moduli spaces for them.
As Andrei Teleman informed me, this problem has been raised by Tyurin.
\end{Rem}
\begin{Loc}
We will say that a pair $(\F,\E)$, consisting of a torsion
free coherent sheaf $\E$ and a saturated non-zero proper subsheaf $\F$,
\it
satisfies the condition
($*$)\rm , if
i) $[\E]\in {\frak F}(\Delta)$, ii) there exists a polarization
$H\in\Delta$ such that
$\alpha$) $\mu_{H}(\F)=\mu_H(\E)$, and $\beta$)
$((\rk\F-1) c_1^2\F - 2\rk\F c_2\F).H^{n-2}\le 0$ and
$((\rk\G-1) c_1^2\G - 2\rk\G c_2\G).H^{n-2}\le 0$, $\G:=\E/\F$.
\end{Loc}
\begin{Lem}
$W^1:=\bigl\{\, x\in (1/r!)\Num(X)\ |\ \exists\ (\F,\E)
\text{ satisfying ($*$) }: x=\xi_{\F,\E}\,\bigr\}$ is a finite set.
\end{Lem}
\begin{pf}
This is an easy adaptation of the proof of Thm.~1.3 in \cite{MW}:
Let $x$ be in $W^1$. Choose a pair $(\F,\E)$ satisfying $(*)$
with $x=\xi_{\F,\E}$.
Define $h:=\max\bigl\{\, (s-1)/2s + (r-s-1)/(2(r-s))\ |\
s=1,...,r\,\bigr\}$, $l:=(r-1)/(2r)$, $k_1:=\max\bigr\{\, c_2.H^{n-2}\
|\ H\in\Delta\,\bigr\}$,  $k_2:=\min
\bigr\{\, c_1^2.H^{n-2}\ |\
H\in\Delta\,\bigr\}$.
Then exactly as in \cite{MW}, p.~105,
one shows that $0\le -x^2.H^{n-2}\le r^2(k_1-lk_2)/(1-h)=:N$.
Observe that  $N$ depends only on $r$,
$c_1$, and
$c_2$. So, it suffices to show that
$$
\bigl\{\, x\in (1/ r!)\Num(X)\ |\ \exists H\in\Delta:\q x.H^{n-1}=0
\q \wedge\q -x^2.H^{n-2}\le N\,\bigr\}
$$
is a finite set.
Again, this can be proved in the same manner as Lemma~1.5 in \cite{MW}.
Indeed, the bilinear form $\langle.\, ,.\rangle_H$ with
$\langle x,y\rangle_H=x.y.H^{n-2}$ depends continuously on $H$, and,
since $H$ is supposed to be a K\"ahler class, it has signature
$(1,\rho(X)-1)$, by the Hodge-Riemann bilinear relations (\cite{GH},
p.~123).
\end{pf}
\subsection{A boundedness result}
The basis of our investigations
is the following
\begin{Prop}
The set
${\frak F}(\Delta)$
is bounded.
\end{Prop}
\begin{pf}
Denote by $W^{1*}$ the set of elements $x\in W^1$ such that
$x.H^{n-1}=0$ for only finitely many polarizations $H\in\Delta$.
For each such $x$, let $w^1(x)$ be the set of $H$ such that
$x.H^{n-1}$ is zero.
We set $w^1:=\bigcup_{x\in W^{1*}} w^1(x)$.
Let $[\E]$ be in ${\frak F}(\Delta)$, such that $\E$ is slope
$H_{\la_0}$-semistable, $\la_0\in\Q$, but fails to fulfill the
assumptions of Proposition~\ref{wichtig1}. Then it is easy to check
that
$H_{\la_0}$ lies in $w^1$.
Let $U_1,...,U_s$ be the connected components of
$\Delta\setminus w^1$. Pick polarizations $A_i\in U_i\cap N^1_{\Q}(X)$,
$i=1,...,s$, and denote by $A_{s+1},...,A_t$ those elements in $w^1$
which are rational.
By
Proposition~\ref{wichtig1}, the concept of slope (semi)stability remains
constant within each $U_i$. So, any $\E$ with $[\E]\in {\frak
F}(\Delta)$
will be slope semistable w.r.t.~one of the polarizations  $A_1,...,A_t$.
\end{pf}
%
\section{Passing through a rational wall}
\subsection{Riemann-Roch}
For any torsion free coherent sheaf $\E$ on $X$, we have its
Chern character $\chh(\E)\in A^*(X)$. We will denote
its homogeneous component of degree $d$ by $\chh_d(\E)$.
We denote by $\td_e$ the degree $e$ part of the Todd character
of the tangent bundle of $X$.
Then, the Riemann-Roch theorem asserts
$$
\chi(\E)\q =\q \sum_{i=0}^n \chh_i(\E).\td_{n-i}.
$$
For any line bundle $\L$ on $X$, we know that
$\chh(\E\otimes\L)=\chh(\E).\chh(\L)$ so that
$$
\chi(\E\otimes\L)\q =\q  {1\over n!}r\L^n+{1\over
(n-1)!}\L^{n-1}.(\chh_1(\E)+r\td_1)+
\cdots+ \chi(\E).
$$
In particular, the Hilbert polynomial
of $\E$ w.r.t.~the ample line bundle $H$ is
$$
P_H(\E)\q =\q
\left({1\over n!}rH^n\right)x^n+\left({1\over
(n-1)!}H^{n-1}.(\chh_1(\E)+r\td_1)\right) x^{n-1}+
\cdots+ \chi(\E).
$$
Define $\ch_d(\E):= \chh_d(\E)+\chh_{d-1}(\E).\td_1+\cdots+ r\td_r$ for
$d=1,...,n$.
To abbreviate notation,
for a subsheaf $\F\subset\E$ and $0\le d\le n$,
we define
$$
\ch_d(\F,\E)\q :=\q {\ch_d(\E)\over \rk\E}-{\ch_d(\F)\over \rk\F}.
$$
\subsection{More walls}
We have already defined a set of walls $w^1$, such that the
concept of slope semistability remains constant between these
walls.
Define $w^2$ as follows: The set of isomorphy classes of sheaves
$\F$ which are saturated subsheaves of sheaves in the family
${\frak F}(\Delta)$, such that
$[(c_1\F/\rk\F)-(c_1/r)].H^{n-1}=0$ for all polarizations in $\Delta$ is
bounded, so that they provide us with a finite set of equations
$\ch_i(\F,\E).H^{n-i}=0$.
We consider only those equations which are non-trivial and let
$w^2$ be set of the respective solutions. Set $w:=w^1\cup w^2$.
By the very definition of $w$, the concept of Gieseker semistability
remains constant within each connected component of $\Delta\setminus w$.
\begin{Rem}
i) The walls in $w\setminus w^1$ do not affect the concept of
slope stability, i.e., the moduli spaces for two polarizations
separated only by a wall in $w\setminus w^1$ will be isomorphic
at least over the open subsets parametrizing slope stable sheaves.
\par
ii)
As we have seen in Example~\ref{zb2},
it is possible that $w$ contains
points which
do not lie in $N^1_{\Q}(X)$. In this case the methods presented in
this section break
down and have only the weak results of Section~\ref{OST}. However, the
reader may check that on some simple manifolds such as $\P_1\times\P_n$,
all the walls are rational. In those cases,
our results completely describe the situation, at least from an
abstract viewpoint. The phenomenon of real walls might explain
the difficulties
encountered by Qin in the definition of walls for higher dimensional
varieties \cite{Q1}.
\end{Rem}
\subsection{The crucial lemma}
\label{Cr}
Suppose that $H_0$ and $H_1$ lie in neighbouring connected components
of $\Delta\setminus w$ which are separated by a \sl rational \rm
polarization
$A$.
We can furthermore assume that there
is an effective $\Q$-divisor $D$ such that $H_1=A+D$ and $H_0=A-D$.
If $X$ is a surface, then in both \cite{MW} and \cite{EG}
the result is based on the fact that there is an integer $l_0$ such
that $\E$ is Gieseker $H_1$- ($H_0$-)(semi)stable if and only if
$\E(l_0D)$ ($\E(-l_0D)$) is Gieseker $A$-(semi)stable.
This result allows one to explore some parameter dependent
(semi)stability concept \sl w.r.t.~the polarization $A$ \rm
such that for different choices of the parameter one obtains
${\mathscr M}_{H_0}$, ${\mathscr M}_{H_1}$, and
${\mathscr M}_{A}$, respectively.
Now, this choice of parameter corresponds in a suitable construction
to the choice of a linearization of a group action. The variation
of the quotients in the latter setting
is well understood. Indeed, this
problem can be appropriately dealt with in the context of master spaces.
In the abstract GIT setting, the construction of master spaces
is carried out in
\cite{Th2}. Examples of master spaces which solve moduli problems
can be found in \cite{OST} and \cite{Sch}.
\begin{Lem}
\label{central}
There is an integer
$l_0$ such that for every $l\ge l_0$ and  every
torsion free coherent sheaf $\E$ with Hilbert form $p$
the following conditions are
equivalent.
\begin{enumerate}
\item $\E$ is Gieseker $H_1$-(semi)stable ($H_0$-(semi)stable).
\item $\E(lD)$ ($\E(-lD)$) is Gieseker $A$-(semi)stable.
\end{enumerate}
\end{Lem}
\begin{pf}
We will explain the proof for $H_1$ in the semistable case.
It is our task to compare the Hilbert polynomials
$P_{H_1}(\E)$ and $P_A(\E(lD))$.
Let $\E$ be a torsion free coherent sheaf
with Hilbert form $p$, and let $\F\subset \E$ be a non-zero proper
subsheaf. One computes
\begin{eqnarray*}
& &\delta(\F,\E,l)(m):= {\chi(\E(lD)\otimes A^m)\over r} -
{\chi(\F(lD)\otimes A^m)
\over\rk\F}\\
 & &\qquad =
B_{n-1} \ch_1(\F,\E).A^{n-1} m^{n-1}
\\
& &\qquad +\left(B_{n-2}^1 A^{n-2}.\ch_2(\F,\E) + B_{n-2}^2 l
A^{n-2}.D.\ch_{1}(\F,\E)\right)m^{n-2}+\\
& &\qquad\qquad\qquad\qquad\vdots\\
& &\qquad +\left(B_{n-i}^1 A^{n-i}.\ch_i(\F,\E)+
\cdots
+
B_{n-i}^i l^{i-1} A^{n-i}.D^{i-1}.\ch_1(\F,\E)\right)m^{n-i}+\\
& &\qquad\qquad\qquad\qquad\vdots\\
& &\qquad +
 B_0^1\ch_n(\F,\E)+\cdots+ B_0^n l^n D^{n-1}.\ch_1(\F,\E).
\end{eqnarray*}
The $B_i^j$ are just some positive constants of no importance.
The coefficient of $m^{n-i}$ in $\delta(\F,\E,l)$ will
be denoted by $\delta_i(\F,\E,l)$.
\par
Assume $\E$ is Gieseker $H_1$-semistable.
First, we know by the $H_1$-semistability of $\E$ and our assumptions on
the walls that $\E$ is at least slope $A$-semistable.
If $\F$ is a non-zero proper subsheaf of $\E$ with
$\ch_1(\F,\E).A^{n-1}>0$, then we see that $\F(lD)$ won't
$A$-desemistabilize $\E(lD)$ for any $l$.
Thus, we can assume that $\ch_1(\F,\E).A^{n-1}=0$.
But the family of all sheaves $\F$ such that there is
a Gieseker $A$-semistable sheaf $\E$ containing $\F$ as a
non-zero proper saturated subsheaf and $\ch_1(\F,\E).A^{n-1}=0$
is bounded. This is important to keep in mind for the
rest of the proof, because it shows that the number
of equations arising in the following is indeed finite, and therefore
one can find an $l_0$ working for all of them.
Now, suppose we have a subsheaf $\F$ of $\E$
such that
$\delta_i(\F,\E,l)=0$ for $i=1,...,j$. By induction we know
that then we must have
$\ch_i(\F,\E).A^{n-i}=0$ for $i=1,...,j$, and
$\ch_i(\F,\E).H_\la^{n-i}=0$ for $i=1,...,j-1$ and any
$H_\la:=A+\la D$ with $\la\in [0,1]$.
If
$\ch_i(\F,\E).H_\la^{n-i}=0$ for any
$H_\la$, $\la\in [0,1]$, then obviously
$A^{n-\iota}.D^{\iota-i}.\ch_i(\F,\E)=0$ for $\iota=i,...,n$.
Therefore, $\delta_{j+1}(\F,\E)=B^1_{n-j-1}
A^{n-j-1}.\ch_{j+1}(\F,\E)+B_{n-j-1}^2\allowbreak l A^{n-j-1}.D.\ch_j(\F,\E).$
If we assume $\ch_j(\F,\E).H_1^{n-j}>0$, then our assumption on
the walls implies that $\ch_j(\F,\E).H_\la^{n-j}>0$ for
all $\la\in (0,1]$. One checks, by choosing $\la$ very small,
that this
forces
$A^{n-j-1}.D.\ch_j(\F,\E)>0$. But then for large $l$,
$\delta_{j+1}(\F,\E,l)>0$, and we don't have to care
about $\F$ any more. If, on the other hand,
$\ch_j(\F,\E).H_1^{n-j}=0$, then our assumption on the walls
shows that $\ch_j(\F,\E).H_\la^{n-j}=0$ for any $\la\in [0,1]$.
The $H$-semistability of $\E$ implies in this case
$\ch_{j+1}(\F,\E).H_1^{n-j-1}\allowbreak\ge 0$. Again using the assumption on
the walls, we will also have $\ch_{j+1}(\F,\E).A^{n-j-1}\ge 0$.
In the present circumstances $\ch_{j+1}(\F,\E).A^{n-j-1} > (=)\ 0$
is equivalent to $\delta_{j+1}(\F,\allowbreak \E,l) > (=)\ 0$.
Either we can stop, or we go on with our induction.
\par
Now, let $\E(lD)$ be $A$-semistable for all $l$ sufficiently large.
First of all, we remark that this implies that $\E$ is
slope $A$-semistable. For any subsheaf $\F\subset\E$ with
$\ch_1(\F,\E).A^{n-1}\allowbreak
>0$, we will also have $\ch_1(\F,\E).H_1^{n-1}>0$.
Hence, only the saturated subsheaves with $\ch_1(\F,\E).A^{n-1}=0$ are
of interest. But these sheaves live again in a bounded family.
Suppose we have a subsheaf $\F\subset\E$ such that
$\ch_i(\F,\E).H_1^{n-i}=0$ for $i=1,...,j-1$ ($j=1$ is allowed). Then,
of course,
$\ch_i(\F,\E).H_\la^{n-i}=0$  for $i=1,...,j-1$  and any $\la\in [0,1]$.
Moreover, $\delta_i(\F,\E,l)=0$ for $i=1,...,j-1$ in this case,
and $\delta_j(\F,\E,l)= B_{n-j}^1 A^{n-j}.\ch_j(\F,\E)$.
Again, $\delta_j(\F,\E,l)>0$ implies $H_1^{n-j}.\ch_j(\F,\E)>0$,
so only the case $\delta_j(\F,\E,l)=0$  matters.
If $j=n$, we get $(\chi(\E)/ r)-(\chi(\F)/\rk\F)=0$, whence
$\F$ does not $H_1$-desemistabilize $\E$.
Otherwise, we look at $\delta_{j+1}(\F,\E,l)=B^1_{n-j-1}
A^{n-j-1}\ch_{j+1}(\F,\E)
+ l B_{n-j-1}^2 A^{n-j-1}.D.\ch_j(\F,\E)$.
If $H_1^{n-j-1}.\ch_j(\F,\allowbreak \E)\allowbreak <0$, then
$H_\la^{n-j-1}.\ch_j(\F,\E)<0$ for all $\la\in (0,1]$.
For small $\la$ this means $A^{n-j-1}.D.\allowbreak \ch_j(\F,\E)<0$.
In this case $\delta_{j+1}(\F,\E,l)<0$ for large $l$,
contradicting our assumptions on $\E$.
\end{pf}
\subsection{Flips between moduli spaces of parabolic sheaves}
As for $\dim X=1$ \cite{Th2}, one can describe the variation of
moduli spaces of parabolic sheaves in terms of GIT flips. Furthermore,
they can be flipped to the corresponding Gieseker moduli space. This
will be worked out in the present section.
\subsubsection*{Parabolic sheaves}
Let $X$ be as before, let $A$ be an ample line bundle on $X$,
and $D\subset X$ an effective divisor.
Fix polynomials $P$, $P_1$,...,$P_k$. Let
$\underline{\alpha}=
(\alpha_{0},...,\alpha_{k})$
be a weight vector with rational entries
$0<\alpha_{0}<\cdots<\alpha_{k}<1$.
A \it parabolic sheaf of weight $\underline{\alpha}$ \rm
is a filtration
$\E=\F_0\supset\F_{1}\supset\cdots\supset\F_{k}\supset
\F_{{k+1}}=\E(-D)$. To shorten notation, we just denote it by
$\E$. Define its \it (parabolic) Hilbert polynomial \rm as
$P_A^{\underline{\alpha}}(\E):= P_A(\E) - \sum_{i=1}^{k+1}
\eps_{i}
P_A(\E/\F_i)$, where $\eps_{i}:=\alpha_{{i}}-
\alpha_{{i-1}}$, $i=1,...,k$, $\alpha_{{k+1}}:=1$.
Given a parabolic sheaf $\E$ of weight $\underline{\alpha}$, every
subsheaf
$\F$ of $\E$ can be viewed as a parabolic sheaf
of weight $\underline{\alpha}$ .
We say that a parabolic sheaf of weight $\underline{\alpha}$ is
\it (semi)stable \rm if for every non-zero proper
subsheaf $\F$ the condition $
P_A^{\underline{\alpha}}(\F)/\rk\F\ (\le\nobreak)\
P_A^{\underline{\alpha}}(\E)/\rk\E$ holds.
Of course, one can also define the \it parabolic slope  \rm
$\mu_A^{\underline{\alpha}}$ of $\E$ and speak of \it slope
semistability\rm .
\par
We restrict our attention to
parabolic sheaves
$\E=\F_0\supset\F_{1}\supset\cdots\supset\F_{k}\supset
\F_{{k+1}}=\E(-D)$
of weight $\underline{\alpha}$ where $P_A(\E)=P$ and $P(\E/\F_{i})
=P_{i}$, $i=1,...,k$.
The moduli space for S-equivalence classes of semistable parabolic
sheaves of weight
$\underline{\alpha}$
was constructed in \cite{MY} and \cite{Y}. Let us denote it by
${\mathscr M}^{\mathop{\rm par}}_A(P,P_1,...,P_k;\underline{\alpha})$.
Below, we will briefly review the construction.
%
\begin{Thm}
\label{Pflips}
Let $P$, $P_1$,...,$P_k$ be as before. Suppose we are given two
weight vectors
$\underline{\alpha}=(\alpha_{0},...,\alpha_k)$
and
$\underline{\alpha}^\p=(\alpha^\p_{0},...,\alpha_{k}^\p)$.
Then ${\mathscr M}_A(P)$, the moduli space of S-equivalence classes of
Gieseker
$A$-semistable torsion free coherent sheaves with Hilbert polynomial
$P$,
${\mathscr M}^{\mathop{\rm par}}_A(P,P_1,...,P_k;\underline{\alpha})$, and
${\mathscr M}^{\mathop{\rm par}}_A(P,P_1,...,P_k;\underline{\alpha}^\p)$
can be all constructed via GIT out of the same quasi-projective scheme,
and, moreover, there is a quasi-projective scheme ${\frak X}$
with an ample line bundle ${\frak L}$ on it and a natural
${\C^*}^{k+1}$-action, and there are linearizations
$\sigma_0$, $\sigma$, and $\sigma^\p$ of this $
{\C^*}^{k+1}$-action in ${\frak L}$ such that
${\frak X}\catqot_{\sigma_0}
{\C^*}^{k+1}
={\mathscr M}_A(P)$, ${\frak X}\catqot_{\sigma} {\C^*}^{k+1}=
{\mathscr M}^{\mathop{\rm par}}_A(P,P_1,...,P_k;\underline{\alpha})$,
and
${\frak X}\catqot_{\sigma^\p} {\C^*}^{k+1}=
{\mathscr M}^{\mathop{\rm par}}_A(P,P_1,...,P_k;\underline{\alpha}^\p)$.
Thus, by the Mum\-ford-Thaddeus principle (\cite{Th2}, \cite{OST},
Part~1), these spaces are related by a sequence of
${\C^*}^{k+1}$-flips.
\end{Thm}
\subsubsection*{Some useful semistability criteria}
Let $W_0,...,W_k$ be finite dimensional $\C$-vector spaces. Define
$W:=W_0\oplus \cdots\oplus W_k$, and let ${\C^*}^k$ act on
$W$ in the following way: The $i$-th factor of ${\C^*}^k$ acts by
scalar multiplication on $W_i$ and trivially on all other summands,
$i=1,...,k$.
In this way, we obtain a linearized action of ${\C^*}^k$
on $\P(W)$. By means of an induction,
one derives the following observation from \cite{OST},
Example~1.2.5.
\begin{Lem}
\label{OST1}
Considering all possible linearizations of the above ${\C^*}^k$-action
on $\P(W)$, one obtains the following polarized quotients
$$
\Bigl( \bigl(\P(W_{\iota_1})\times\cdots\times
\P(W_{\iota_\kappa})\bigr),
\bigl[\O(a_1,...,a_\kappa)\bigr]\Bigr).
$$
Here, $\{\, \iota_1,...,\iota_\kappa\,\}$ can be any subset
of $\{\, 0,...,k\,\}$, and $(a_1,...,a_\kappa)$ any tuple of
positive integers.
\end{Lem}
Consider a reductive algebraic group $G$ and 
representations $\rho_i\colon G\allowbreak\lra\allowbreak\GL(W_i)$,
$i=0,...,k$. The direct sum of these representations defines an
$\O_{\P(W)}(1)$-linearized action of
$G$ on $\P(W)$. We also have $\O_{\P(W_i)}(1)$-linearized actions
of $G$ on $\P(W_i)$, $i=0,...,k$, and for a point
$[v_i]\in\P(W_i)$ and a one parameter subgroup $\la\colon\C^*\lra G$
we let $\mu_i([v_i],\la)$ be minus the weight of the
induced $\C^*$-action on the fibre of $\O_{\P(W_i)}(1)$ over the
point $\lim_{z\lra \infty}\la(z)\cdot [v_i]$.
\begin{Prop}
\label{OST2}
Let $w=[v_0,...,v_k]\in\P(W)$ be a point, and let $(\nu_1,...,\nu_\mu)$
be the indices with $v_{\nu_j}\neq 0$, $j=1,...,\mu$.
Then the following conditions are equivalent:
\begin{enumerate}
\item $w$ is $G$-semistable w.r.t.~given linearization.
\item There exist non-negative integers $l_{\nu_1},...,l_{\nu_\mu}$,
not all zero, such that
for any one parameter subgroup $\la\colon
\C^*\lra G$
$$
l_{\nu_1}\mu_{\nu_1}([v_{\nu_1}], \la)+\cdots+
l_{\nu_\mu}\mu_{\nu_\mu}([v_{\nu_\mu}], \la)\q\ge\q 0.
$$
\end{enumerate}
\end{Prop}
\begin{Rem}
In view of Lemma~\ref{OST1}, the second conditions means that we find a
linearization of the ${\C^*}^k$-action such that the image of
$w$ in the corresponding polarized quotient is $G$-semistable
w.r.t.~the induced linearization.
\end{Rem}
\begin{pf} We observe that the hypothesis that $G$ have no characters
in Section~1.2.\ of \cite{OST} only assures that the linearization of
$G$ is unique. In the proofs, this assumptions is never used.
So, we can apply \cite{OST}, Thm.~1.4.1, to
prove the assertion by induction. The details are left to the reader.
\end{pf}
\subsubsection*{A "baby" master space construction}
In this section, we explain the proof of Thm.~\ref{Pflips}.
To avoid excessive indices and formulas, we will only treat the case
$k=0$ which is the only one we will need for our applications. Using
the semistability criteria given above, the reader will have no
difficulty to extend the proof to the case of arbitrary $k$.
We need to fix a Poincar\'e sheaf ${\frak P}$ on $\Pic X\times X$.
\par
First of all, we may choose an integer $m_0$ such that for every
$m\ge m_0$ and every torsion free coherent sheaf $\E$ which is
either slope $A$-semistable or which appears in a parabolic sheaf
of either weight $\underline{\alpha}$ or $\underline{\alpha}^\p$
\begin{itemize}
\item $H^i(X, \E(mA))=0$ for $i=1,...,n$.
\item $\E(mA)$ is generated by global sections.
\item The same holds for $\E_{|D}(mA)$.
\end{itemize}
Moreover, let ${\frak A}\subset\Pic X$ be the union of all components
      containing elements of the form $[\det\E]$.
\begin{itemize}
\item
      Then $\L(r mA)$
      is globally generated and without higher cohomology for every
      $[\L]\in {\frak A}$.
\end{itemize}
As usual, we consider the Quot scheme ${\frak F}$
of equivalence classes
of
quotients $q\colon V\otimes \O_X(-mA)\lra \E$ where $\E$ is a coherent
$\O_X$-module with Hilbert polynomial $P$.
Furthermore, there is a
universal flag
$$
V\otimes \pi_X^*\O_X(-mA)\lra {\frak E}_{\frak F}\lra {\frak E}_{{\frak
F}|{\frak F}\times D}
$$
over
${\frak F}\times X$.
Let $U_0$ be the set of points $\bigl[q\colon
V\otimes\O_X(-mA)\lra\E\lra \E_{|D}\bigr]$ for which $\E$ is Gieseker
$A$-semistable, let $U_{\underline{\alpha}}$ and
$U_{\underline{\alpha}^\p}$ be
the sets for which $\E\supset\E(-D)$ is a semistable
parabolic sheaf of weight $\underline{\alpha}$ and
$\underline{\alpha}^\p$,
resp., and $U:=U_0\cup U_{\underline{\alpha}}\cup
U_{\underline{\alpha}^\p}$. The sheaf $\pi_{{\frak F}*}\bigl({\frak
E}_{\frak F}\otimes\pi_X^*\O_X(mA)\bigr)$ is locally free of rank
$P(m)$,
and the sheaf
$\pi_{{\frak F}*}\bigl({\frak
E}_{{\frak F}|{\frak F}\times D}\otimes\pi_X^*\O_X(mA)\bigr)$
is locally free of rank, say, $R$.
The scheme $U$ can the be mapped $\SL(V)$-equivariantly to
$$
\P\Bigl(\underline{\Hom}\bigl( \bigwedge^rV\otimes \O_{\frak A},
\pi_{{\frak A}*}\bigl({\frak P}\otimes
\pi_X^*\O_X(mA)\bigr)^\vee\Bigr)\times
\P\Bigr( \bigwedge^R\bigl(V\otimes H^0(\O_X(mA))\bigr)^\vee\Bigr).
$$
Let $\P_{\frak A}$ be the first factor of this product,
and $\P_R$ the second.
Choose some ample sheaf ${\frak H}_{\frak A}$ on ${\frak A}$,
so that ${\frak L}_{\frak A}:= \O_{\P_{\frak A}}(1)\otimes \pi_{\frak
A}^*
{\frak H}_{\frak A}$  is ample. The sheaf $\pi_{\P_{\frak A}}^*
{\frak L}_{\frak A}^{\otimes a}\otimes\pi_{\P_R}^*\O_{\P_R}(b)$
on $\P_{\frak A}\times \P_R$ will be denoted by $\O(a,b)$.
Denote by $U^\p_0$ the set of $\SL(V)$-semistable points
w.r.t.~the linearization in $\O(1,0)$. Then $U_0$ is mapped injectively
and properly to $U^\p_0$, and for suitable choices
of $(a,b)$ and $(a^\p,b^\p)$, the sets $U_{\underline{\alpha}}$
and $U_{\underline{\alpha}^\p}$ get immersed into the sets
$U^\p_{\underline{\alpha}}$
and $U^\p_{\underline{\alpha}^\p}$ of points which are
$\SL(V)$-semistable
w.r.t.~the linearization in $\O(a,b)$ and $\O(a^\p,b^\p)$, respectively.
Altogether, we obtain an injective and proper map of $U$ to
$U^\p:=U_0^\p\cup U^\p_{\underline{\alpha}}\cup
U^\p_{\underline{\alpha}^\p}
$.
It is now clear that the moduli spaces, we are interested in are
obtained from $U$ by dividing out $\SL(V)$ for different linearizations.
To understand the assertion about the $\C^*$-flips, we proceed as
follows. Define ${\frak R}$ as the projective bundle
over ${\frak Q}$ associated to the vector bundle
$$
\pi_{{\frak F}*}\underline{\Hom}\Bigl(\det\bigl({\frak E}_{\frak
F}\otimes\pi_X^*\O_X(mA)\bigr), (\det\times{\id}_X)^*{\frak P}\Bigr)
\oplus \bigwedge^R \bigl(V\otimes H^0(\O_X(mA))\bigr)\otimes\O_{\frak
F},$$
$\det\colon {\frak F}\lra{\frak A}$ being associated with the
family ${\frak E}_{\frak F}$, and
${\frak S}$ the projective bundle
$$
\P\Bigl(
\underline{\Hom}\bigl( \bigwedge^rV\otimes \O_{\frak A},
\pi_{{\frak A}*}({\frak P}\otimes
\pi_X^*\O_X(mA))\bigr)^\vee
\oplus
\bigwedge^R\bigl(V\otimes H^0(\O_X(mA))\bigr)^\vee\otimes\O_{\frak A}
\Bigr)
$$
over ${\frak A}$. One has the natural morphism ${\frak t}\colon
{\frak R}\lra {\frak S}$ (compare \cite{OST}, Section~2.4). There are
natural
$(\SL(V)\times\C^*)$-actions
on ${\frak R}$ and ${\frak S}$, and ${\frak t}$ is equivariant.
The $\SL(V)$-action is canonically linearized, and we can choose
linearizations $s_0$, $s$, and $s_1$ of the $\C^*$-action such
that the polarized quotients are
\begin{eqnarray*}
{\frak S}\catqot_{s_0}\C^* &=&\bigl(\P_{\frak A}, [{\frak L}_{\frak
A}]\bigr);\\
{\frak S}\catqot_{s}\C^* &=&\bigl(\P_{\frak A}\times
\P_R,[\O(a,b)]\bigr);\\
{\frak S}\catqot_{s^\p}\C^* &=&\bigl(\P_{\frak A}\times
\P_R,[\O(a^\p,b^\p)]\bigr).
\end{eqnarray*}
Let $U^{\p\p\p}_0$, $U^{\p\p\p}_{\underline{\alpha}}$, and $
U^{\p\p\p}_{\underline{\alpha}^\p}$ be the respective sets of
$(\SL(V)\times\C^*)$-semistable points, and let $U^{\p\p\p}$ be their
union.
Their preimages $U_0^{\p\p}$,
$U^{\p\p}_{\underline{\alpha}}$, and $
U^{\p\p}_{\underline{\alpha}^\p}$ under ${\frak t}$ coincide with
the preimages of
$U_0$,
$U_{\underline{\alpha}}$, and $
U_{\underline{\alpha}^\p}$ under the bundle map ${\frak R}\lra {\frak
F}$. Thus, the union $U^{\p\p}$ of these sets maps finitely to
$U^{\p\p\p}$. By general properties of good quotients, the quotient
${\frak Y}:= U^{\p\p\p}\catqot \SL(V)$ is an open subset of the
projective scheme
${\frak S}\catqot \SL(V)$, and ${\frak X}:= U^{\p\p}\catqot \SL(V)$
maps finitely to ${\frak Y}$; call the corresponding map ${\frak z}$.
Both, ${\frak X}$ and ${\frak Y}$ inherit $\C^*$-actions, and
${\frak z}$ is equivariant w.r.t.~them. By construction
and the "commutation principle"
(e.g., \cite{OST}, Sect.~1.3.1),
the $\C^*$-action
on ${\frak Y}$ is linearized in an ample
line bundle ${\frak L}_{\frak Y}$ such that suitable manipulations
of this linearization will yield ${\frak
S}\catqot_{s_0}(\SL(V)\times\C^*)$ and so on as quotients.
Pulling back these linearizations to ${\frak X}$ gives us ${\frak L}$,
$\sigma_0$,
$\sigma$, and $\sigma^\p$ as asserted.\hfill\qed
\subsection{The proof of the Main Theorem}
We return to the setting of Section~\ref{Cr} and choose
some $l$ for which Lemma~\ref{central} holds. For a torsion
free coherent sheaf $\E$ and $\beta\in [0,1]$, we set $P_A^\beta(\E):=
(1-\beta) P_A(\E(-lD)) + \beta P_A(\E(lD))$, and call $\E$
\it $\beta$-(semi)stable\rm, if and only if $P_A^\beta(\F)/\rk\F\ (\le)\
P_A^\beta(\E)/\rk\E$ for any non-trivial proper subsheaf $\F$. In
Lemma~\ref{central}, we have
seen that a torsion free coherent sheaf $\E$ with Hilbert form $p$ is
$H_1$-($H_0$-)(semi)stable if and only if $\E$ is
$1$-($0$-)(semi)stable. But as the proof of Lemma~\ref{central} shows,
we can choose $\beta_1$ close to one and $\beta_2$ close to zero,
so that we will also have that $\E$ is
$H_1$-($H_0$-)(semi)stable if and only if $\E$ is
$\beta_1$-($\beta_0$-)(semi)stable. As a corollary to the existence of
moduli
of parabolic bundles (the r\^ole of $\E$ is the last section will now
be played by $\E(lD)$ and that of $D$ by $2lD$), for any
$\beta\in
(0,1)$, there exists a projective moduli scheme
${\mathscr M}^\beta_A(p)$ of S-equivalence classes of $\beta$-semistable
torsion free coherent sheaves with Hilbert form $p$, and as we have seen
in~\ref{central}
${\mathscr M}_A^{\beta_i}(p)\cong {\mathscr M}_{H_i}(p)$, for $i=0,1$.
Therefore, the main theorem is a direct consequence of
Theorem~\ref{Pflips}.\hfill\qed
\section{Passing through an arbitrary wall}
\label{OST}
Let $H_0$ and $H_1$ be two polarizations, and ${\frak F}(H_0)$ and
${\frak F}(H_1)$ be the set of isomorphy classes of torsion free
coherent sheaves which are slope $H_0$-semistable and
slope $H_1$-semistable, respectively. Let $H$ be an arbitrary
polarization and write $\O_X(m)$ for $\O_X(mH)$.
Since both
${\frak F}(H_0)$ and ${\frak F}(H_1)$
are bounded, we can find a complex vector spaces
$V$ and an integer $m_0$ such that any sheaf $\E$ whose isomorphy
class belongs to either ${\frak F}(H_0)$ or ${\frak F}(H_1)$
can be embedded into $V\otimes \O_X(m)$ for all $m\ge m_0$.
We denote by ${\frak Q}$ the Quot scheme of all submodules
of $V\otimes \O_X(m_0)$ with Hilbert form $p$. Strictly speaking, this
is a fine moduli space of $\delta$-stable pairs $(\E,\phi)$, $\phi\in
\Hom(\E,V\otimes\O_X(m_0))$, for some large polynomial $\delta$.
But as its universal property shows, it is isomorphic to a Quot scheme
and, in particular, does not depend on the choice of a polarization.
\par
Fix a Poincar\'e sheaf ${\frak P}$ on $\Pic X\times X$, and
let ${\mathscr M}_{H_i/{\frak P}/V\otimes\O_X(m_0)}(p)$ be the master space
of S-equivalence classes of semistable ${\frak P}$-oriented pairs
$(\E,\eps,\phi)$ \cite{OST} where $\E$ is a torsion free coherent sheaf
with
Hilbert polynomial $P_{H_i}(n)=p(H_i^{\otimes n})$, for all $n\in\N$,
$\eps\colon\det\E\lra
{\frak P}_{|\{[\det\E]\}\times X}$ is a homomorphism, and
$\phi\in\Hom(\E,V\otimes\O_X(m_0))$, $i=1,2$.
As proved in \cite{OST}, there are natural $\C^*$-actions on these
master spaces.
Suitably linearized, these $\C^*$-actions give rise to sequences of
$\C^*$-flips which begin with a fibration
$\pi_i\colon {\frak M}_{i}\lra {\mathscr
M}_{H_i}(p)$ and end in ${\frak Q}$.
The fibre of $\pi_i\colon {\frak M}_{i}\lra {\mathscr
M}_{H_i}(p)$ over the ismorphy class of a stable sheaf $\E$
is just $\P\bigl(\Hom(\E, V\otimes\O_X(m_0))^\vee\bigr)$.
Therefore, we have
shown that the fibrations
$\pi_0\colon {\frak M}_{0}\lra {\mathscr
M}_{H_0}(p)$ and
$\pi_1\colon {\frak M}_{1}\lra {\mathscr
M}_{H_1}(p)$ can be created by means of $\C^*$-flips out of the
Quot scheme ${\frak Q}$.
%

\par
\bigskip
Bar-Ilan University\par
Department of Mathematics and Computer Science\par
Ramat-Gan, 52900\par
Israel\par
\medskip
and\par
\medskip
Universit\"at GH Essen\par
FB 6 Mathematik und Informatik\par
D-45117 Essen\par
Deutschland\par
\smallskip
\it E-mail address: \tt schmitt1@@cs.biu.ac.il
\end{document}